%% file: main.tex
\documentclass[letterpaper, 10 pt, conference]{ieeeconf}  

\IEEEoverridecommandlockouts                              

\usepackage{comment}
\usepackage{cite}
\usepackage{svg}
\usepackage{tikz}
\usepackage{amsmath,amssymb,amsfonts}
\usepackage{algorithmic}
\usepackage{graphicx}
\usepackage[hyperfootnotes=false]{hyperref}
\usepackage{textcomp}
\usepackage{subcaption}
\usepackage{siunitx}
\usepackage{bm}
\usepackage{physics}
\usepackage[enable]{easy-todo}
\usepackage{url}
\usepackage{acro}

\input{definitions}

\usepackage{lipsum}

\makeatletter
\newcommand{\blankfootnote}[1]{%
  \bgroup
  \renewcommand{\thefootnote}{}%
  \footnotetext{\refstepcounter{footnote}#1}%
  \egroup
}
\makeatother

\begin{document}
\title{Expected Time-Optimal Control: a Particle Model Predictive Control-based Approach via Sequential Convex Programming}
\author{Kazuya Echigo, Abhishek Cauligi, and Beh\c{c}et A\c{c}{\i}kme\c{s}e \IEEEmembership{Fellow, IEEE}
\thanks{The work was supported by Air Force Office of Scientific Research grant FA9550-20-1-0053. The research was carried out in part at the Jet Propulsion Laboratory, California Institute of Technology, under a contract with the National Aeronautics and Space Administration (80NM0018D0004).}
\thanks{Kazuya Echigo and Beh\c{c}et A\c{c}{\i}kme\c{s}e are with the William E. Boeing Department of Aeronautics and Astronautics, University of Washington, Seattle, WA 98195, USA. (e-mail: kazuyae@uw.edu). Abhishek Cauligi is with the Jet Propulsion Laboratory, California Institute of Technology.}
}

\newcommand{\spann}{\mathrm{span}}
\newcommand{\ourproblem}{OURPROBLEM}
\newcommand{\nameofourproblem}{NAMEOFOURPROBLEM}
\newcounter{pcnt}
\newtheorem{thm}{Theorem}[section]
\newtheorem{lem}{Lemma}[section]
\newtheorem{asn}{Assumption}[section]
\newcommand{\Problem}{\stepcounter{pcnt}\nid{{\sf \large \textbf{Problem
    \arabic{pcnt}.}} }}
\pagestyle{empty} 
\maketitle
\thispagestyle{empty} 

\begin{abstract}
In this paper, we consider the problem of minimum-time optimal control for a dynamical system with initial state uncertainties and propose a sequential convex programming (SCP) solution framework. We seek to minimize the expected terminal (mission) time, which is an essential capability for planetary exploration missions where ground rovers have to carry out scientific tasks efficiently within the mission timelines in uncertain environments. 
Our main contribution is to convert the underlying stochastic optimal control problem into a deterministic, {\em numerically tractable}, optimal control problem.
To this end, the proposed solution framework combines two strategies from previous methods: i) a partial model predictive control with consensus horizon approach and ii) a sum-of-norm cost, a temporally strictly increasing weighted-norm, promoting minimum-time trajectories. Our contribution is to adopt these formulations into an SCP solution framework and obtain a numerically tractable stochastic control algorithm. We then demonstrate the resulting control method in multiple applications: i) a closed-loop linear system as a representative result (a spacecraft double integrator model), ii)  an open-loop linear system (the same model), and then iii) a nonlinear system (Dubin's car). 
\end{abstract}

\input{introduction}

\input{problemformulation2}
\input{consensuscontrol}

\input{proposedidea}
\input{numerical}
\input{conclusion}

\textcolor{black}{
\bibliographystyle{IEEEtran}
\bibliography{reference}
}
\end{document}

%% file: definitions.tex
\DeclareAcronym{MPC}{
  short=MPC,
  long=model predictive control,
}


%% file: introduction.tex
\section{Introduction}
This paper considers an optimal control problem of minimizing the expected terminal (mission) time of a deterministic dynamical system from an uncertain initial state to a prescribed destination within a sequential convex programming (SCP) framework. The primary motivation for this problem is the upcoming CADRE Lunar rover mission~\cite{DeLaCroixRossiEtAl2024}, in which the rover has to reach a destination as quickly as possible in order to carry out scientific measurements effectively. The problem of deterministic minimum-time optimal control has been studied extensively across aerospace applications, including spacecraft rendezvous and proximity operations~\cite{Harris20142, Eren2017, Richard2019}.

\begin{figure}[t!]
        \centering   
        \includegraphics[width=0.43\textwidth]{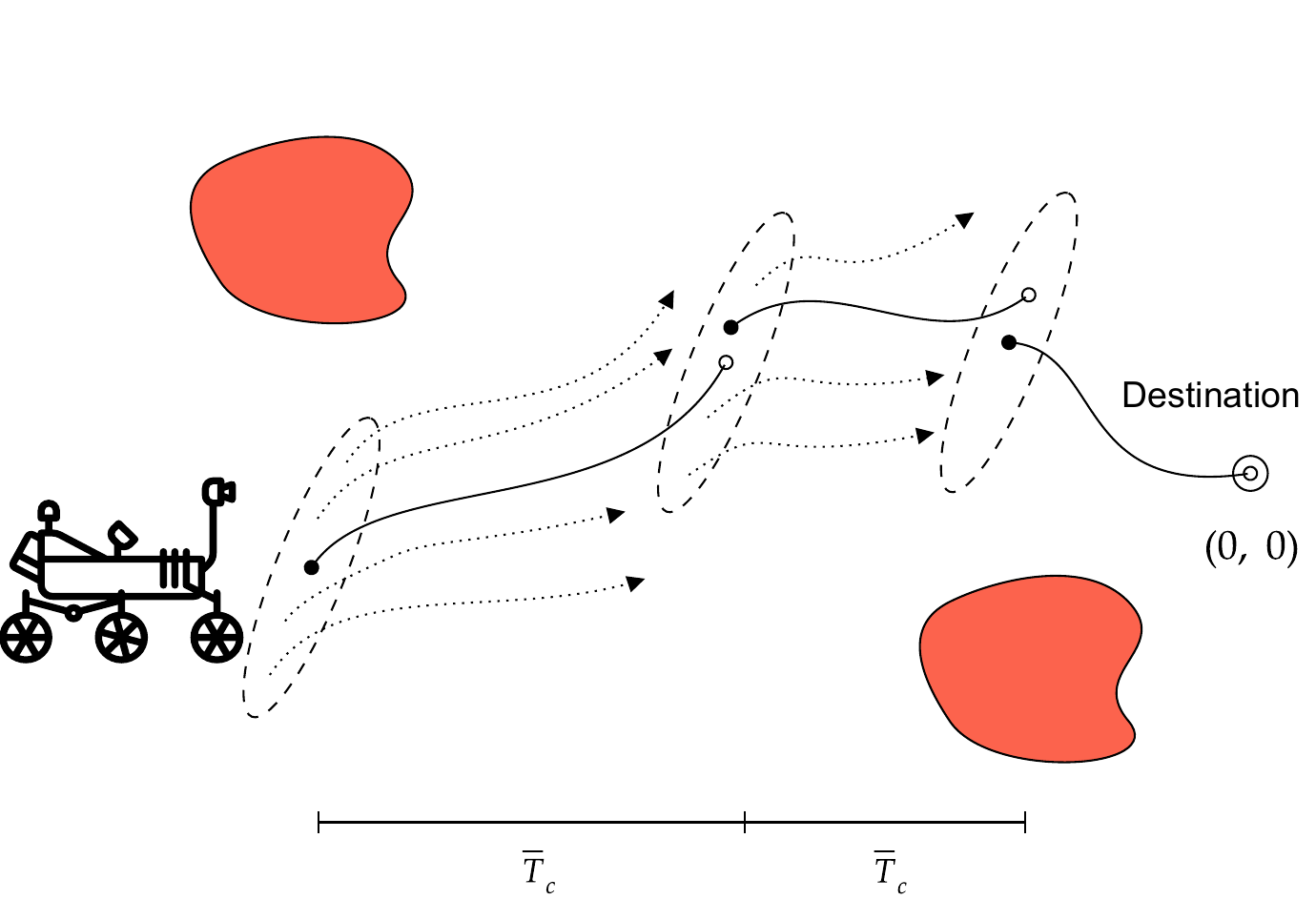}
    \caption{
Overview: The proposed framework minimizes the expected terminal (mission) time to reach a prescribed destination under uncertainties at the initial state. Solving it in a closed loop enables to consider measurement uncertainties.}
\vspace{-1em}
\label{fig:enter-label}
\end{figure}

In such planetary surface rover missions, finding planning policies to minimize the terminal time under uncertainties is crucial. Moreover, autonomous rovers for such missions must be equipped with the ability to replan efficiently online in order to handle multiple sources of uncertainty that may impede carrying out scientific operations, including surface dynamics and rover intrinsic parameters. Due to the lack of a reliable global localization capability for upcoming Lunar rover mission concepts~\cite{CauligiSwanEtAl2023,DaftryEtAl2023}, a major source of uncertainty is that of measurement uncertainty in the rover's initial position estimate. Motivated by this shortcoming, we provide a robust framework to generate a path with minimum expected terminal time under measurement uncertainties. As described in Figure~\ref{fig:enter-label}, solving the aforementioned problem in a closed loop can reasonably consider them. 

At present, there are several research works to tackle a minimum expected time optimal control problem, which originated from~\cite{Heath1987}. They include deriving analytical feedback control laws using the Hamilton-Jacobi-Bellman equation~\cite{Anderson2013} and solving the expected shortest path problem in discrete state and action spaces with control uncertainty~\cite {GUILLOT2020148}. However, the former is limited to a certain specific system, and the latter approach cannot handle state uncertainty and is prone to the curse of dimensionality. Incorporating state uncertainty into an optimal control framework while maintaining algorithmic flexibility remains an open question.

In light of the preceding discussion, we propose an efficient approach that minimizes the expected time to reach the prescribed destination under uncertainty. We reformulate the proposed stochastic optimal control problem as a deterministic and numerically tractable optimization problem by introducing a particle~\ac{MPC} approach with a consensus control horizon~\cite{DyroHarrisonEtAl2021} and a sum-of-norm cost~\cite{BakoChenEtAl2011}. In the end, we apply SCP to solve the resulting deterministic optimization problem. SCP, an iterative optimization framework capable of problems with a general set of nonconvexities, is ideal for fast and efficient onboard trajectory planning~\cite{Malyuta2022}. 

Particle~\ac{MPC} maintains the flexibility offered by SCP in stochastic systems, which is capable of handling arbitrary probability distributions~\cite{DyroHarrisonEtAl2021, BlackmoreOnoEtAl2010, HOMEMDEMELLO2014}. In contrast to earlier studies that focused on either sharing control inputs solely at the initial time step of each scenario (resulting in a too-optimistic solution) or throughout all time steps (resulting in an infeasible solution for each scenario to reach the prescribed destination), the authors in~\cite{DyroHarrisonEtAl2021} introduce a promising flexible approach with a tunable consensus horizon. Denoting the consensus control horizon as $\overline{T}_c$, they achieve the best control policies for the next $\overline{T}_c$ steps. In this paper, we employ that idea in a closed-loop fashion by aligning the consensus horizon with a feedback control period. It should be noted that other methods to approximate state uncertainties by stochastic distributions, such as Gaussian distributions~\cite{ZhuAlonsoMora2019} and Gaussian Mixture Models~\cite{Weissel2009}, are not suitable for the aforementioned problem as they also enforce state uncertainties to share control inputs throughout all time steps.

The flexibility offered by the particle~\ac{MPC} approach further opens the door to extending the sum-of-norm cost for a stochastic system, which was originally provided as the relaxation of deterministic minimum-time optimal control problem~\cite{BakoChenEtAl2011}. In deterministic systems, the minimum-time optimal control problem eventually requires executing computationally intensive bisection searches~\cite{SinghSingla2007}. The authors of\cite{BakoChenEtAl2011} approximate that problem in linear deterministic systems as a sparse optimal control problem and further relaxed it by the sum-of-norm cost. It has been proven that the sum-of-norm cost promotes a sparsity in many control research fields~\cite{candes2008, Nagahara2016, Kazu2023LCSS}. The major advantage of that method is that the resulting cost function is convex, which can be solved efficiently using an off-the-shelf convex solver with the guarantee of the global optimal convergence~\cite{DomahidiChuEtAl2013}. 

{\em Statement of Contributions:} The primary contribution of this paper is the framework capable of more efficiently evaluating the expected terminal time for the aforementioned stochastic system with the flexibility of handling various nonlinear stochastic constraints. To efficiently solve the problem, we present a method based on combining the particle~\ac{MPC} framework and the sum-of-norm cost. Specifically, the terminal time of each sampled scenario is approximated with a sum-of-norm cost, then the mean of the sum-of-norm cost functions is calculated. This new cost function maintains convexity, making it suitable for SCP. Next, we provide proof that the mean of the sum-of-norm cost serves as a reasonable convex relaxation of a sparse optimal control problem even in the proposed framework. Finally, we demonstrate our proposed approach can outperform the deterministic minimum-time optimal control through numerical simulations for linear and nonlinear systems. 

The remainder of this paper is organized as follows: Section II outlines a method for approximating the proposed stochastic optimal control problem as a deterministic nonlinear optimization problem, introducing the particle~\ac{MPC} approach. Section III then introduces the sum-of-norm cost to relax the minimum expected time optimal control problem. Lastly, Section IV presents some numerical analysis to compare with benchmark approaches for deterministic systems.

%% file: problemformulation2.tex
\section{Problem Formulation and Particle Model Predictive Control Approach}
In this section, we review how to reformulate a stochastic optimal control problem as a deterministic one using a particle~\ac{MPC} approach with a consensus control horizon. For nonlinear dynamics with uncertainties at the initial time, we address the problem of minimizing the expected time to transfer from an initial state to a prescribed final destination. In this study, we apply a receding horizon control strategy to tackle the optimal control problem, optimizing control inputs in several future time steps. Repeating this strategy in a closed-loop fashion can take an updated state estimate and corresponding uncertainties into account. Furthermore, to reasonably handle those uncertainties appearing in the problem formulation, we apply a particle-based representation of state uncertainty together with a parameterized consensus horizon, as in~\cite{DyroHarrisonEtAl2021}.

We begin by defining the nonlinear dynamics as follows:
\begin{equation}
x(k+1) = f(x(k), u(k)), \nonumber
\end{equation}
where $x(t) \in \mathbb{R}^{n_x}$ is the system state and $u(t) \in \mathbb{R}^{n_u}$ is a deterministic control input. The dynamics $f: \mathbb{R}^{n_{x}} \times \mathbb{R}^{n_{u}} \times \mathbb{R}^{n_p} \rightarrow \mathbb{R}^{n_{x}}$ are assumed to be a continuously differentiable function. Further, without a loss of generality, we assume that the function has an equilibrium point at $(x,u) = 0$. 

With the dynamics described, we now formulate a stochastic control problem for reaching a prescribed destination from an initial uncertain state, as follows: 
\newline
\textbf{Problem 1}:
\label{pro:problem1}
\begin{align}
    \underset{u}{\textrm{minimize}} \,\, & \mathbb{E}(t_f) \nonumber \\
    \textrm{subject to} 
    \,\,& x(k+1) = f(x(k), u(k))\nonumber \\ 
    \,\,&  x(1) \sim \mathcal{X}_1 \nonumber\\
    \,\,& x(t_f) = 0,\nonumber
\end{align}
where we assume that the state vector at the initial time follows a discrete distribution $\mathcal{X}_1$. 

%% file: consensuscontrol.tex

We approximate the uncertainty distribution in the system's initial states with randomly sampled particles, each representing an independent and identically distributed realization of $\mathcal{X}_1$ and a potential path of the system starting from that realization. This approach enables us to consider the mutual effects of dynamics and uncertainty following arbitrary probability distributions. 
We recover a deterministic control policy for this problem by enforcing consensus on the control variables across each of the particles. While the traditional methods can either consider a complete consensus or a one-step consensus, the authors of~\cite{DyroHarrisonEtAl2021} introduced a concept of a tunable consensus horizon, enabling a calibration between those two extreme strategies. Denoting the number of samples as $m$, we introduce that method to yield the following sampled optimization problem:
\newline
\textbf{Problem 2}:
\label{pro:problem2}
\begin{align*}
    \underset{u^i(k),\overline{u}(k), T_i}{\textrm{minimize}} \,\, &\mathbb{E} (T_i),\\
    \textrm{subject to} 
    \,\,& x^i(k+1) = f(x^i(k), u^i(k)),\,\, i \in \mathbb{N}_1^m, \,\, k \in \mathbb{N}_{1}^{\infty}, \\
    \,\,&  x^i(1)  = x^i_1, \,\, i \in \mathbb{N}_1^m, \\
    \,\,&  x^i(T_i)  = 0, \,\, i \in \mathbb{N}_1^m, \\
    \,\,&  u^i(k) \in \mathcal{U}, \,\, i \in \mathbb{N}_1^m, \,\, k \in \mathbb{N}_1^{T_i - 1},\\
    \,\,&  u^i(k) = \overline{u}(k), \,\, i \in \mathbb{N}_1^m, \,\, k \in \mathbb{N}_{1}^{\overline{T_c} - 1}, 
\end{align*}
where $x^i(1)$ is an initial state randomly sampled from the discrete distribution $\mathcal{X}_1$ and $\overline{T}_c$ represents a consensus time horizon with corresponding control $\overline{u}(k)$. The set of integers from $a$ to $b$ is denoted by $\mathbb{N}_a^b$.

%% file: proposedidea.tex
\section{Sum-of-Norm Approximation}
In this section, we outline our approach to solve the minimum expected time optimal control problem as previously introduced. 
The terminal time in each scenario is approximated as a sparse cost function of the state trajectory and then further relaxed as a sum-of-norm cost function. 

We begin by assuming a control horizon $\Gamma$, which exceeds the optimized terminal time per each sample, meaning $\Gamma > T_i^{*}, \,\, \forall i = 1, ..., m.$ Additionally, to simplify the discussion, we further assume that the consensus time horizon $\overline{T_c} <  T_i^{*}$ for any $i$. We represent $u^{*}$ and $x^{*}$ as the optimal trajectories for control and state, respectively, in Problem~\hyperref[pro:problem2]{2}. They are collections of control and state trajectories tailored for individual sampled scenarios, with the minimum expected terminal time to reach the destination. For this solution, it is reasonable to assume that once the vehicle reaches the destination for each scenario, the state remains there. Therefore, the solution state trajectory $x^{*}$ should satisfy the condition as follows:
\begin{eqnarray}
 x^i(k)^{*} \begin{cases}
 \neq 0, \quad (i \in \mathbb{N}_1^m, \,\, k \in \mathbb{N}_1^{T_i^{*} - 1}), \\
 =  0, \quad  (i \in \mathbb{N}_1^m, \,\, k \in \mathbb{N}_{T_i^{*}}^{\Gamma}),  
 \end{cases}
\label{eq:stateprofilecondition}
\end{eqnarray}
where $x^i(k)^{*}$ denotes the state at a specific discrete time $k$ within the $i^\textrm{th}$ sampled trajectory of $x^{*}$. Under this assumption, the optimized terminal time $T_i^{*}$ defined for each scenario is equal to the $L_0$ norm of the following vector:
\begin{eqnarray}
&&[\|x^i(1)^{*}\|_2, ..., \|x^i(\Gamma)^{*}\|_2] \nonumber \\
&& = [\|x^i(1)^{*}\|_2, ..., \|x^i(T_i^{*})^{*}\|_2, 0_{1\times(\Gamma - T_i^{*})}], \nonumber
\end{eqnarray}
which is a $\Gamma$ dimensional vector that contains the elements of the $L_2$ norms of the state at each time step. At the optimum, since $x^{*}$ satisfied the aforementioned assumption, all elements from the first to the $T_i^{*}$ th positions of the vector are non-zero, while all elements after $T_i^{*}$ th positions of the vector are zero. Therefore, minimizing the terminal time for this scenario is identical to minimizing the number of non-zero elements in this vector, i.e., yielding the sparsest solution. Together with the fact that $\mathbb{E}(T_i)$ in Problem~\hyperref[pro:problem2]{2} is represented as $\frac{1}{m} \sum_{i = 1}^{m}T_i$, introducing it equivalently transforms Problem~\hyperref[pro:problem2]{2} into the following optimization problem: 
\newline
\textbf{Problem 3}:
\label{pro:problem3}
\begin{align*}
    \underset{u^i(k),\overline{u}(k), T_i}{\textrm{minimize}} \,\, & \frac{1}{m}\sum_{i = 1}^{m} \| [\|x^i(1)\|_2, ..., \|x^i(\Gamma)\|_2] \|_0,\\
    \textrm{subject to}
    \,\,& x^i(k+1) = f(x^i(k), u^i(k)),\,\, i \in \mathbb{N}_1^m, \,\, k \in \mathbb{N}_{1}^{\Gamma}, \\
    \,\,&  x^i(1)  = x^i_1, \,\, i \in \mathbb{N}_1^m, \\
    \,\,&  u^i(k) \in \mathcal{U}, \,\, i \in \mathbb{N}_1^m, \,\, k \in \mathbb{N}_1^{T_i - 1},\\
    \,\,&  u^i(k) = \overline{u}(k), \,\, i \in \mathbb{N}_1^m, \,\, k \in \mathbb{N}_{1}^{\overline{T_c} - 1}.
\end{align*}
The non-convexity of the $L_0$ norm, characterized by its discontinuity, generally yields an NP-hard optimization problem~\cite{Natarajan1995}. To tackle this, it is well known that the $L_1$ norm cost promotes sparsity and works well as its convex relaxation~\cite{Candes2005}. Given that, in this paper, we introduce its reasonable convex approximation as in~\cite{BakoChenEtAl2011}, which gives:
\newline
\textbf{Problem 4}:
\label{pro:problem4}
\begin{align*}
\underset{u^i(k),\overline{u}(k), T_i}{\textrm{minimize}} \,\, & \frac{1}{m}\sum_{i = 1}^{m} \sum_{k = 1}^{\Gamma} \omega(k) \|x^{i}(k)\|_2 \\
    \textrm{subject to}
    \,\,& x^i(k+1) = f(x^i(k), u^i(k)),\,\, i \in \mathbb{N}_1^m, \,\, k \in \mathbb{N}_1^{T_i}, \\
    \,\,&  x^i(1)  = x^i_1, \,\, i \in \mathbb{N}_1^m, \\
    \,\,&  u^i(k) \in \mathcal{U}, \,\, i \in \mathbb{N}_1^m, \,\, k \in \mathbb{N}_1^{T_i - 1},\\
    \,\,&  u^i(k) = \overline{u}(k), \,\, i \in \mathbb{N}_1^m, \,\, k \in \mathbb{N}_{1}^{\overline{T_c} - 1},
\end{align*}
where $\omega(k) \in \mathbb{R}^{+}$ is a weight parameter that monotonically increases with respect to $k$. Note that those weight parameters are usually introduced to enhance the sparsity of the optimal solution~\cite{candes2008}. 
The cost function introduced by this approximation is convex since $\omega(k)$ is regarded as a predefined parameter within the optimization framework, allowing the formulation to be embedded into an iterative convex subproblem solver such as SCP.

With the approximated optimization problem stated, we now introduce the technical result in Theorem 3.1, which proves that the optimal solution to Problem~\hyperref[pro:problem4]{4} satisfies \eqref{eq:stateprofilecondition} for any $\Gamma$. That is, it not only proves that the $L_1$ norm relaxation satisfies the assumption to transform Problem~\hyperref[pro:problem2]{2} into Problem~\hyperref[pro:problem4]{4}, 
but also implies that the optimal solution to Problem~\hyperref[pro:problem4]{4} does not depend on the control horizon $\Gamma$ that is usually defined as a parameter before computation. 
\begin{thm}
    Assume that the optimal trajectories of Problem~\hyperref[pro:problem4]{4} are given as $u^{**}$ and $x^{**}$, respectively. Then, for any $i'$, when $x^{i'}(k)^{**}$ becomes $0$ at $k = k^{i'}_0, \,\, 1 \leq k^{i'}_0 \leq \Gamma$, $ x^{i'}(k)^{**}$ remains $0$ at any $k$ such that $ k^{i'}_0 \leq k \leq \Gamma$.
\end{thm}
\begin{proof}
If there exist feasible state profiles for any sample $i'$, such that it is equal to $x^{i'}(k)^{**}$ till $k = k^{i'}_0$ and remains at $0$ after $k = k^{i'}_0$, $x^{i'}(k)^{**}$ should be equal to $0$ after $k = k^{i'}_0$ so that it has the optimal cost. To show that this is indeed the situation, we begin by assuming there exists the following control trajectory, $u^{1}$:
\begin{eqnarray}
    u^i(\cdot)^{1} = \begin{cases}
    u^i(\cdot)^{**} &(i \neq i')\\
    \begin{cases}
    u^{i'}(k)^{**} &\quad (1 \leq k \leq k^{i'}_0 - 1)\\ 
    0           &\quad (k^{i'}_0 \leq k).
    \end{cases}
    &(i = i')
    \end{cases}
    \nonumber
\end{eqnarray}
Then we denote $x^{1}$ as a corresponding trajectory propagated by $u^{1}$. Since $u^{i'}(\cdot)^{1}$ is taking the same control as $u^{i'}(\cdot)^{**}$ till $k = k^{i'}_0 - 1$ that drives the states to the $0$, $x^{i'}(\cdot)^{1}$ is equal to $0$ at $k = k^{i'}_0$. Given that, the aforementioned assumption that $f(x,u)$ has an equilibrium point at $(x,u) = (0, 0)$ yields 
\begin{eqnarray}
    x^i(\cdot)^{1} = 
    \begin{cases}
    x^i(\cdot)^{**} &(i \neq i')\\
    \begin{cases}
    x^{i'}(k)^{**} &\quad (1 \leq k \leq k^{i'}_0 - 1)\\ 
    0           &\quad (k^{i'}_0 \leq k).
    \end{cases} &(i = i')\\
    \end{cases} \label{eq:definofstate}
\end{eqnarray}
Now, we shall analyze the cost function provided by that. Given that $u^{**}$ is the optimal solution trajectory of Problem~\hyperref[pro:problem4]{4}, the cost functions given by $u^{**}$ and $u^{1}$ have
\begin{gather}
\frac{1}{m}\sum_{i = 1}^{m} \sum_{k = 1}^{\Gamma} \omega(k) \|x^{i}(k)^{**}\|_2  \leq \frac{1}{m}\sum_{i = 1}^{m} \sum_{k = 1}^{\Gamma} \omega(k) \|x^{i}(k)^{1}\|_2 \nonumber \\
\Leftrightarrow \frac{1}{m} \sum_{k = 1}^{\Gamma} \omega(k) \|x^{i'}(k)^{**}\|_2  \leq \frac{1}{m} \sum_{k = 1}^{\Gamma} \omega(k) \|x^{i'}(k)^{1}\|_2 \nonumber \\
\Leftrightarrow \frac{1}{m} \sum_{k = k^{i'}_0 + 1}^{\Gamma} \omega(k) \|x^{i'}(k)^{**}\|_2  \leq \frac{1}{m} \sum_{k = k^{i'}_0 + 1}^{\Gamma} \omega(k) \|x^{i'}(k)^{1}\|_2 \nonumber\\
\leq 0. \nonumber            
\end{gather}
The second inequality follows from the definitions of $x^{1}$ as in \eqref{eq:definofstate}. Furthermore, $x^{i'}(k)^{1}$ is taking the same value as $x^{i'}(k)^{**}$ on $[1 \,\, k^{i'}_0 - 1]$ and then remaining to $0$ on $[k^{i'}_0\,\, \Gamma]$, which gives the third and fourth inequality, respectively. Consequently, the assumption that $\omega(k) > 0$ implies $x^{i'}(k)^{**} = 0$ on $[k^{i'}_0\,\, \Gamma]$. 
\end{proof}

To this end, the proposed cost function is convex, being suitable for SCP-based algorithms. It should be noted that Problem~\hyperref[pro:problem4]{4} only requires linearization of the dynamics; in practice, this structure has been observed to lead to strong convergence properties when using SCP. 
It should be noted that the size of variables in Problem~\hyperref[pro:problem4]{4} scales linearly with respect to $m$. Therefore, addressing all propagated trajectories of multiple scenarios through a one-shot optimization problem might seem computationally expensive. However, we argue that the computational feasibility of these sampling-based approaches is increasing significantly. This is primarily due to the recent development of variance reduction techniques like importance sampling~\cite{HOMEMDEMELLO2014} and GPU acceleration methods~\cite{ChariKamathEtAl2024}.

%% file: numerical.tex
\begin{figure}[t!]
  \centering
\includegraphics[width=0.42\textwidth]{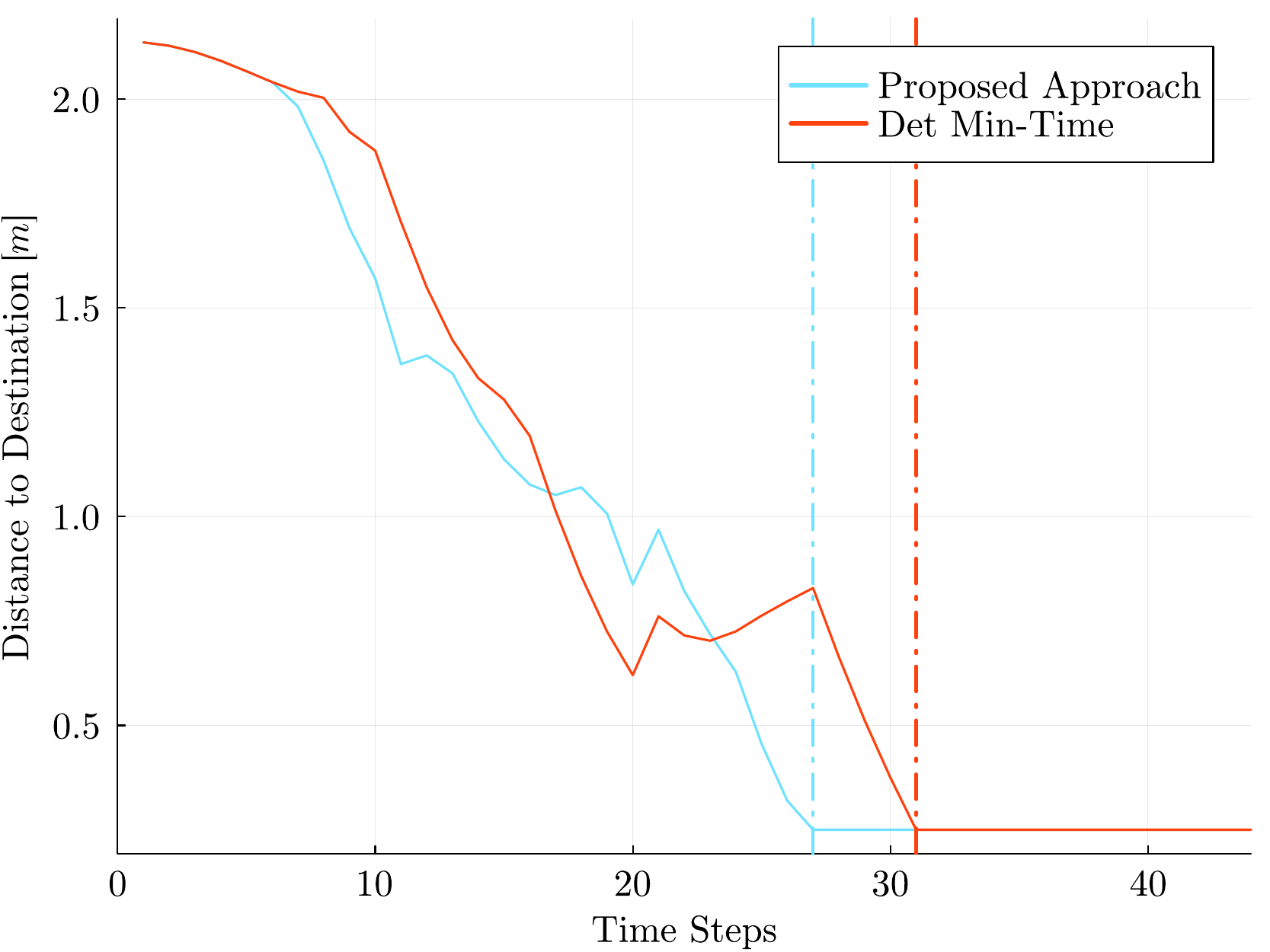}
    \caption{
    The $L_2$ norm distance time histories to the prescribed destination using the proposed approach and the deterministic minimum-time optimal control. It also displays the convergence time steps of those two methods (blue and red dotted lines). This indicates the capability of our proposed approach to achieve faster convergence in contrast to the deterministic controller.}
    \label{fig:representative_example}
    \vspace*{-0.3cm}
\end{figure} 

\section{Numerical Examples}  

In this section, we demonstrate the efficacy of our proposed approach through simulation results on two benchmark problems. After showing a representative result in a double integrator model used for spacecraft dynamics in~\cite{ MoteEgerstedtEtAl2020,CauligiChakrabartyEtAl2022}, the performance of this approach is then demonstrated statistically in Monte Carlo analyses for that double integrator model and a nonlinear Dubin's car model. Our approach is mainly compared against a popular alternative: minimizing the terminal time in deterministic systems as in~\cite{BakoChenEtAl2011}.

Within those Monte Carlo analyses, we apply the proposed approach in an open-loop fashion to facilitate numerical analysis: that is, we implement the proposed approach for the system where uncertainty exists at the initial point. We then propagate randomly sampled scenarios taken from $\mathcal{X}_1$ using the consensus control given by the method. When a discrete time step reaches $\overline{T}_c$, we apply the approximated deterministic minimum-time optimal control for each propagated point. The method from~\cite{BakoChenEtAl2011} is used to compute that deterministic minimum-time optimal control for the double integrator model. For the Dubins car model, we apply the Dubins path~\cite{Dubin1957} for simplicity. On the other hand, in comparing the proposed approach with the deterministic minimum-time optimal control, we apply that deterministic approach in lieu of the proposed approach under the assumption that the deterministic initial point is given as the mean of $\mathcal{X}_1$.

We solve the optimization problems using the SCP toolbox~\cite{Malyuta2022} and ECOS~\cite{DomahidiChuEtAl2013} as the underlying convex optimization problem solver in the Julia programming language. 

\subsection{Representative Result (Closed Loop)}

In this subsection, we demonstrate our approach on a spacecraft double integrator model using two feedback controllers: the proposed approach using Problem~\hyperref[pro:problem4]{4} and the deterministic minimum-time optimal controller with an approximated cost $ \sum_{k = 1}^{\Gamma} \omega(k)\|x(k)\|_2$ as in~\cite{BakoChenEtAl2011}. Note that the frequency of feedback control is set to $1 / \overline{T_c}$. The dynamics of the system, state and control constraints, and necessary parameters for them are summarized as follows:
\begin{gather}
    x(k+1) = \begin{bmatrix}
        I & T_s I \\
        0 & I
    \end{bmatrix} x(k) + \begin{bmatrix}
        \frac{1}{2}T_s^2 I \\
        T_s I
    \end{bmatrix} u(k), \,\, \| u(t) \|_{\infty}\leq 1, \nonumber \\
    T_s = \frac{8}{\Gamma -1}, \,\, \Gamma = 60, \,\, \overline{T}_c = 8, \,\, \omega(k) = k, \nonumber \\
    \mathcal{X}_1 = \mathcal{N}(x_1, \Sigma_1), 
    \,\, x_1 = [2, 1, 0, 0]^\top, \,\, \Sigma_1 = \begin{bmatrix}
I& 0_{2\times2}\\
0_{2\times2} & I
\end{bmatrix} \nonumber \\
m = 30,  \,\, x_f = [0, 0, 0, 0]^\top, \nonumber 
\end{gather}
where $x \in \mathbb{R}^{4}$, $u \in \mathbb{R}^{2}$, and $I$ is denoted as a $2$ dimensional identity matrix. Furthermore, every time we apply a feedback control, Gaussian measurement uncertainty is added to a true propagated state, whose mean $\mu$ and its covariance matrix $\Sigma$ are given as follows: 
\begin{gather*} \mu= 0_{4\times1}, \,\, 
\Sigma = \begin{bmatrix}
0.1 I& 0_{2\times2}\\
0_{2\times2} & 0.1 I
\end{bmatrix}.
\end{gather*}
Therefore, denoting the propagated true state as $x_\textrm{pro}$ gives the updated $\mathcal{X}_1$ in Problem~\hyperref[pro:problem1]{1} of the form:
\begin{gather}
    \mathcal{X}_1 = \mathcal{N}(\mu + x_\textrm{pro}, \Sigma). \nonumber
\end{gather}
Finally, once the distance to $x_f$ becomes less than $0.25$, it is regarded as having converged. Figure~\ref{fig:representative_example} shows the result of one representative scenario, showing that the proposed approach converges faster than the deterministic minimum-time optimal controller. 

\begin{figure}[t!]
        \centering    
    \includegraphics[width=0.45\textwidth]{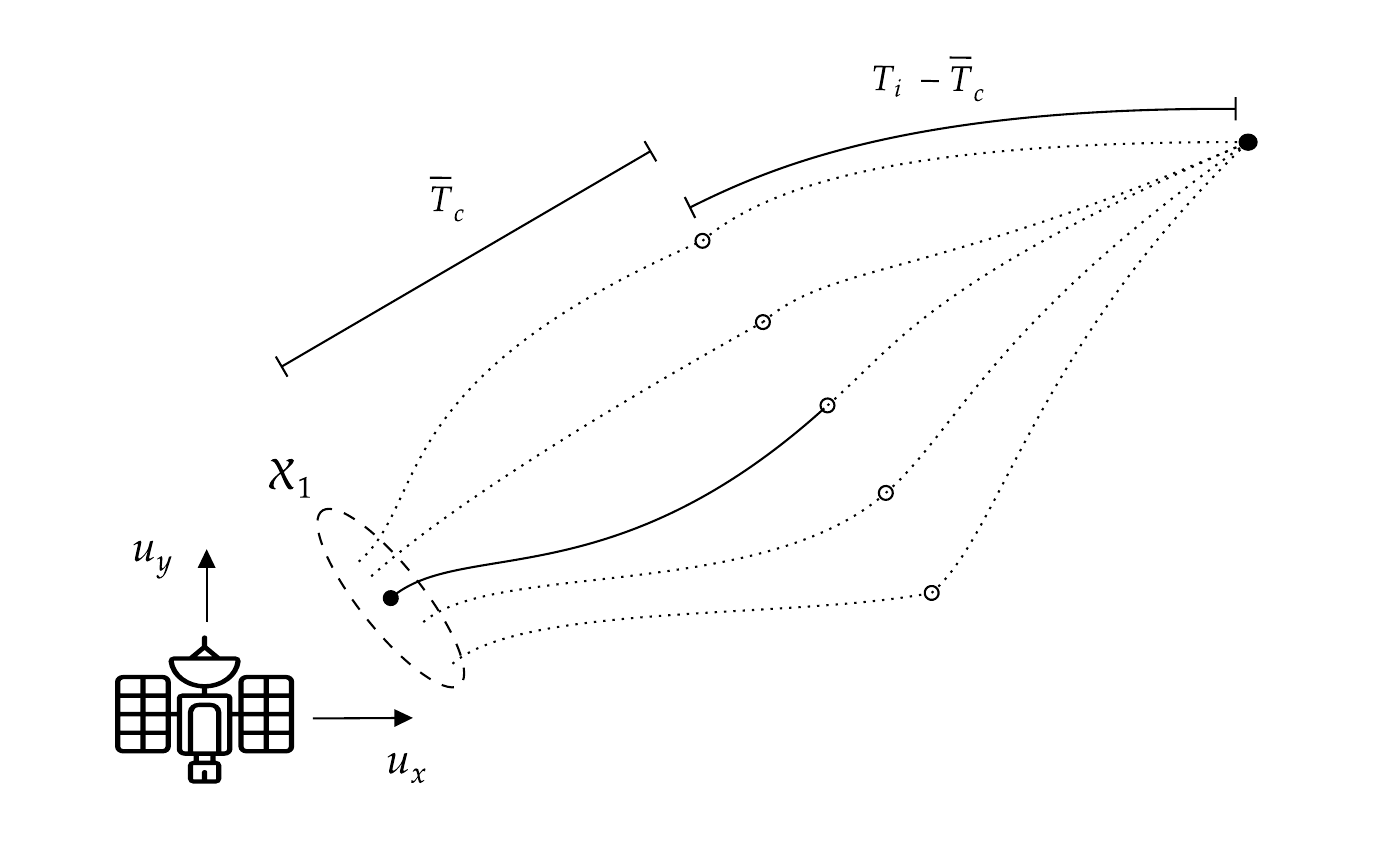}
    \caption{4D free-flying spacecraft system.}
\label{fig:SC_open}
\end{figure}

\begin{figure}[t!]
  \centering
\includegraphics[width=0.48\textwidth]{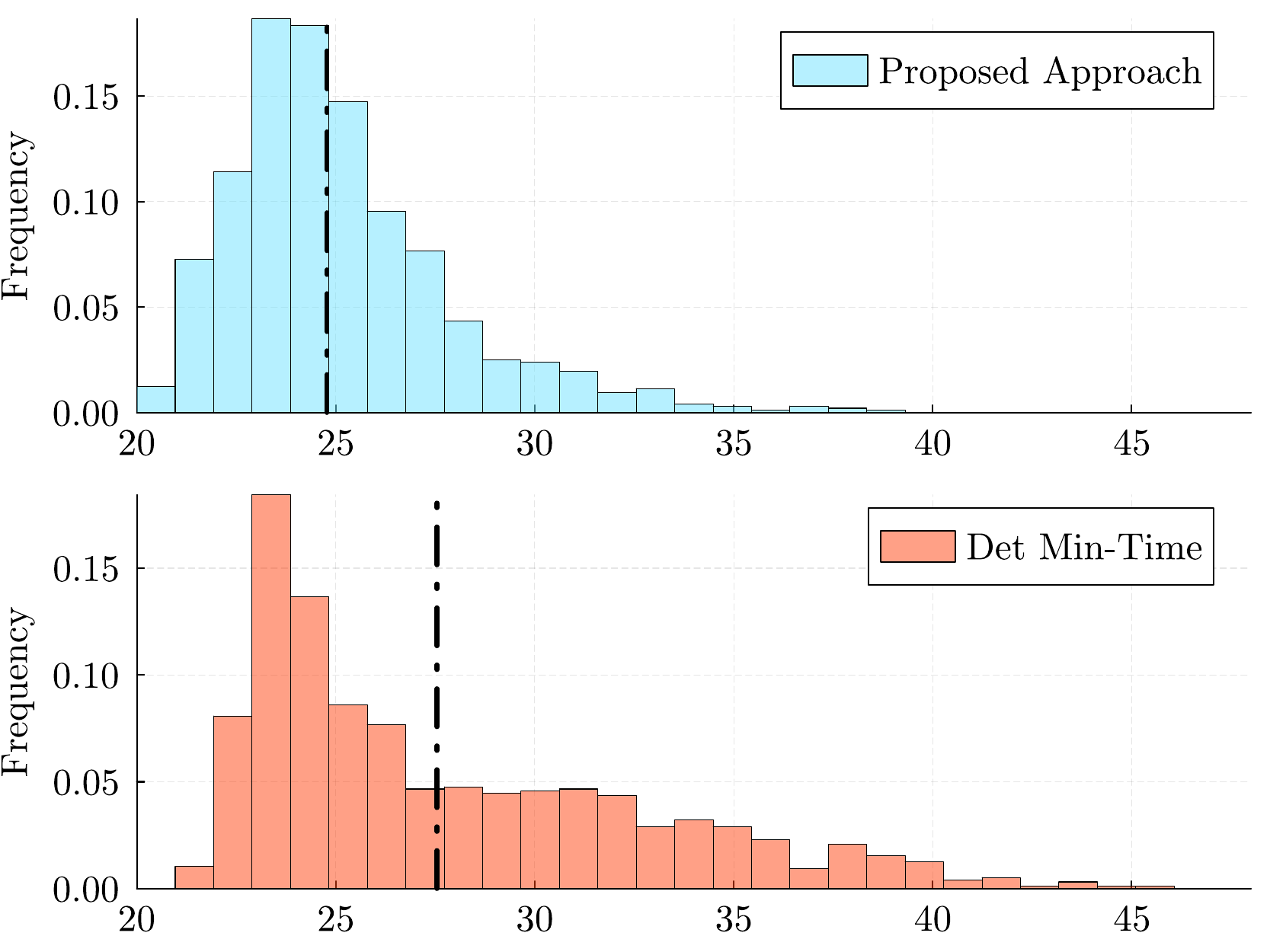}
    \caption{Resulting distributions of terminal time-steps for each trajectory from the proposed approach and the deterministic minimum-time optimal control. Each figure also displays the mean terminal time-steps (black dotted line). The mean terminal time-steps for the proposed approach is less than that for the other.}
    \label{fig:MCresults_simple}
\vspace{-0.2cm}
\end{figure}

\subsection{Linear Systems (Open Loop)}
The capability of our approach is now presented in detail with the aforementioned deterministic minimum-time optimal controller. As depicted in Figure~\ref{fig:SC_open}, those methods are tested in an open-loop fashion for the sake of clarity. Therefore, uncertainties only appear at the beginning. Based on the system presented in the last subsection, the following parameters are updated: 
\begin{gather}
T_s = \frac{6}{\Gamma - 1}, \,\, \Gamma = 40, \,\,  \overline{T}_c = 10, \nonumber  \\
    \mathcal{X}_1 = \mathcal{N}(x_1, \Sigma_1), \,\, \Sigma_1 = \begin{bmatrix}
0.2 I& 0_{2\times2}\\
0_{2\times2} & 0.2 I
\end{bmatrix}, \nonumber
\end{gather}  
Given that, we perform a Monte Carlo analysis with $1000$ samples. Note that a convergence threshold is set to $\SI{1e-3}{}$.

Figure~\ref{fig:MCresults_simple} shows the histograms of the distribution of terminal time-steps. The mean terminal time-steps for the proposed approach is $24.8$, whereas that for the deterministic approach is $27.5$. As expected, our proposed approach outperforms the other when minimizing the expected value of terminal times under uncertainty. 

In~\cite{BakoChenEtAl2011}, it has previously been observed that the sufficient condition of weights to recover the original deterministic minimum-time optimal control was highly conservative, and all of the monotonically increasing weight sequences provided in the paper achieved the minimum-time optimal solution. Given that, Table~\ref{tab:weightparams} describes the result of a further sensitivity analysis of weight parameters for our stochastic system, which suggests that the aforementioned flexibility also holds even after the proposed stochastic extension. 


\begin{table}[t]
    \centering
    \begin{tabular}{lc}
    \hline
   \textbf{Weights} &  \textbf{Mean Terminal Time-Steps} \\ \hline
    $\omega_\textrm{const} := 1$ & $24.8$\\
    $\omega_\textrm{lin} := k$ & $24.8$\\
    $\omega_\textrm{log} := \log_e(k)$ & $24.8$\\
    $\omega_\textrm{quad} := k^2$ & $24.8$\\
   \hline
    \end{tabular} 
    \caption{Weight parameter comparison. The mean terminal time-steps remain the same even when changing the weight parameter.}
    \label{tab:weightparams}
\end{table}  

\begin{table}[t!]
    \centering
    \begin{tabular}{cccccccccc} 
    \hline
    \textbf{Consensus Control Horizon} & 2& 6 & 10 & 14 & 18 \\ \hline
    \textbf{Mean Terminal Time-Steps} & 22.4 & 22.6 & 24.8& 28.3 & 32.6 \\
    \hline
    \end{tabular} 
    \caption{Sensitivity analysis on consensus control horizon. We see that longer $T_c$ results in a more conservative solution, whereas shorter $T_c$ allows for actions suited to each sampled scenario, leading to an overly optimistic solution.}
    \label{tab:sensitivity_analysis}
\end{table}

\begin{table}[t!]
    \centering
    \begin{tabular}{ccccc}
    \hline   
   \textbf{Parameter} & $\omega_\textrm{vc}$ & $\omega_\textrm{tr}$ & $\delta_\textrm{tol}$ & $q$ \\ \hline
   \textbf{Value} & $\SI{1e2}{}$ &$\SI{1e-2}{}$ & $\SI{5e-3}{}$ & Inf \\ \hline
    \end{tabular} 
    \caption{Algorithm parameters for PTR}
    \label{tab:PTRparams}
    \vspace*{-0.2cm}
\end{table}

Finally, Table~\ref{tab:sensitivity_analysis} shows another sensitivity analysis of consensus control horizon $\overline{T_c}$. 
As expected, we see that increasing $\overline{T_c}$ imbues additional conservatism into the controller and results in an expected mean final time.
Thus, a practitioner can tune $\overline{T_c}$ to extract the desired trade-off between conservatism and performance for the controller.

\subsection{Nonlinear Systems (Open Loop)}

Now, we focus on minimizing the expected time to get to a prescribed destination with a Dubins car model under initial state uncertainty. Since the dynamics are nonlinear, we employ SCP to solve Problem~\hyperref[pro:problem4]{4}. Specifically, we use the penalized trust region (PTR) SCP algorithm~\cite{taylor2021,szmuk2019}, whose corresponding parameters are given in Table~\ref{tab:PTRparams}. It should be noted that the proposed approach is not limited to the PTR; it is also applicable to any other SCP-based algorithms~\cite{Malyuta2022}.

\begin{figure}[t!]
  \centering
\includegraphics[width=0.48\textwidth]{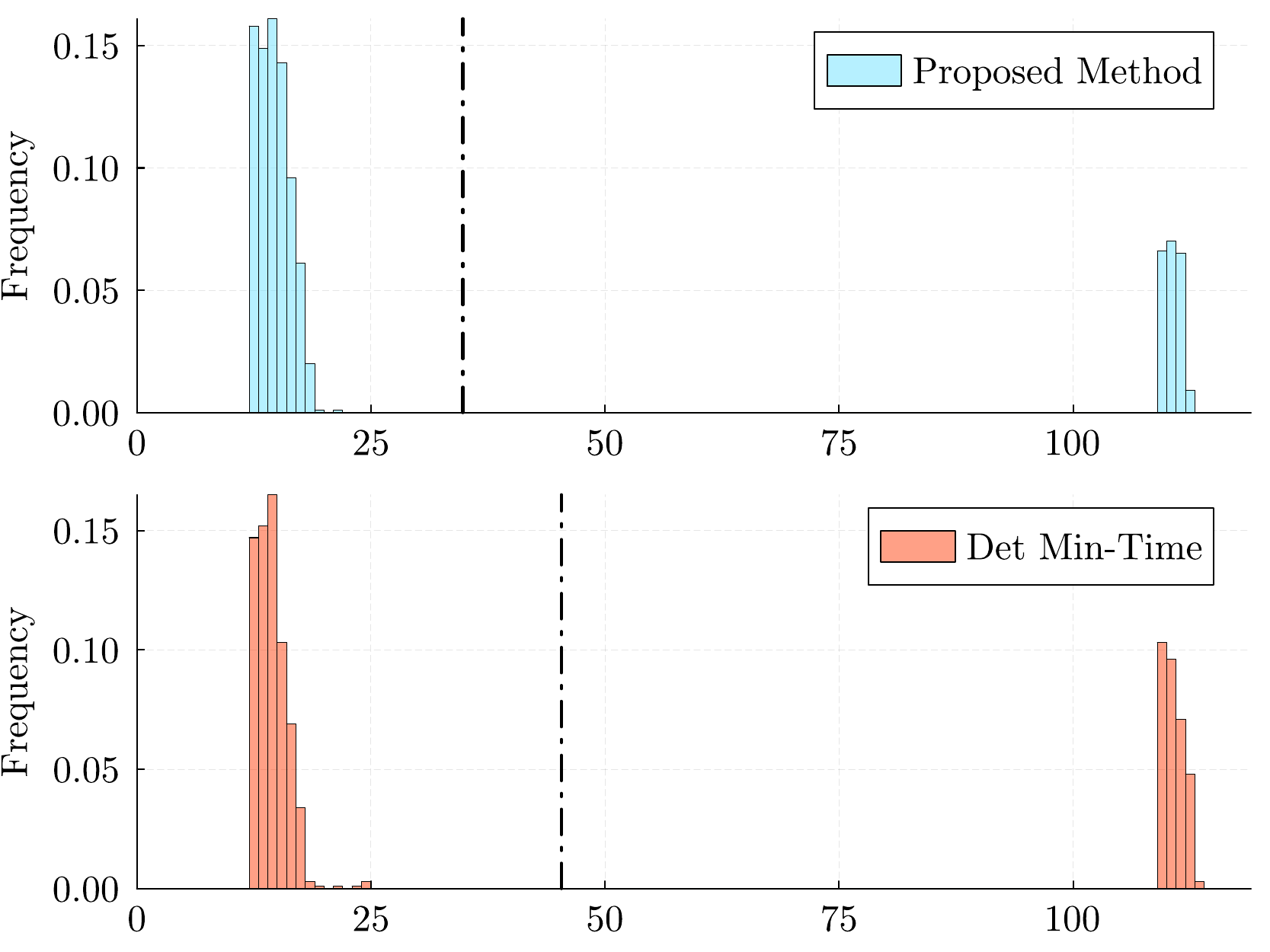}
    \caption{
    Distributions of the terminal time-steps for each trajectory show that our proposed approach outperforms the deterministic minimum-time optimal control approach in terms of the mean terminal time-step (black dotted line).
    }
\label{fig:MCresults_dubin}
\end{figure} 

\begin{table}[t]
    \centering
    \begin{tabular}{ccc}
    \hline   
    &  \multicolumn{2}{c}{\textbf{Samples / (Mean Terminal Time-Steps)}}   \\ 
   \cline{2-3}  
   \textbf{Methods} &  Time Steps $< 20$ & Time Steps $ \geq 20$  \\ \hline
   \textbf{Proposed Approach} & 790 / (14.4) & 210 / (112.6)\\
   \textbf{Det Min-Time} & 679 / (14.4) & 321 / (112.8)\\ \hline
    \end{tabular} 
    \caption{The summary of Monte Carlo simulations}
    \label{tab:monte_dubin}
\vspace*{-0.2cm}
\end{table}

Inspired by several aerospace rover projects~\cite{DeLaCroixRossiEtAl2024}, a three-degree-of-freedom (3-DOF) state space model is considered of the form:
\begin{gather}
x(k+1) = x(k) + T_s * [v\sin(\theta), v\cos(\theta), \phi]^\top, \nonumber
\end{gather}
where $x = [r_x, r_y, \theta] \in \mathbb{R}^3$ is the vehicle's state and $u = [v, \phi] \in \mathbb{R}^2$ is the vehicle's control input. Control constraints, the number of samples, and the other necessary parameters are defined as follows:
\begin{gather}
     0 \leq v(t) \leq 0.5, \,\, |\phi(t)| \leq 5, \,\, T_s = \frac{0.4}{\Gamma -1}, \nonumber \\ 
     x_f = 0_{3\times1}, \,\, \Gamma = 32, \,\,
    \overline{T}_c = 12, \,\, \omega(k) = k, \,\, m = 20. \nonumber 
\end{gather}
Initial state uncertainty is described as follows:
\begin{gather}
    \mathcal{X}_1 = \mathcal{N}(x_1, \Sigma_1), \,\, x_1 = [0, -1, 0]^\top, \,\, \Sigma_1 = \textrm{diag}([0, 0.1, 0]). \nonumber
\end{gather}
In the PTR algorithm, we normalized all variables in advance to a range between $0$ and $1$. The PTR algorithm also requires an initial trajectory guess, for which a linear interpolation is chosen from $x_1$ to $x_f$ with a constant control vector $u(t) = [0.25, 0]^\top$ for every $t$. Then, the resulting control profiles are evaluated via $\SI{1000}{}$ Monte Carlo simulations, whose convergence criterion is set to $\SI{1e-3}{}$.

\begin{figure}[t!]
     \centering
     \begin{subfigure}[b]{0.235\textwidth}
         \centering         \includegraphics[width=\textwidth]{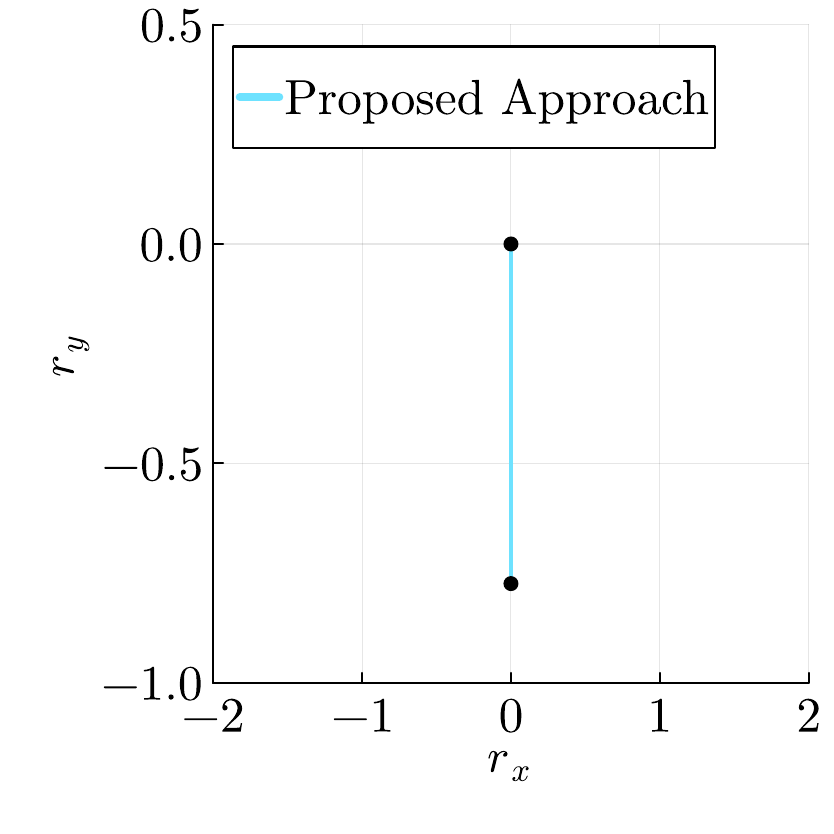}
         \label{fig:traj}
     \end{subfigure}
     \begin{subfigure}[b]{0.235\textwidth}
         \centering
        \includegraphics[width=\textwidth]{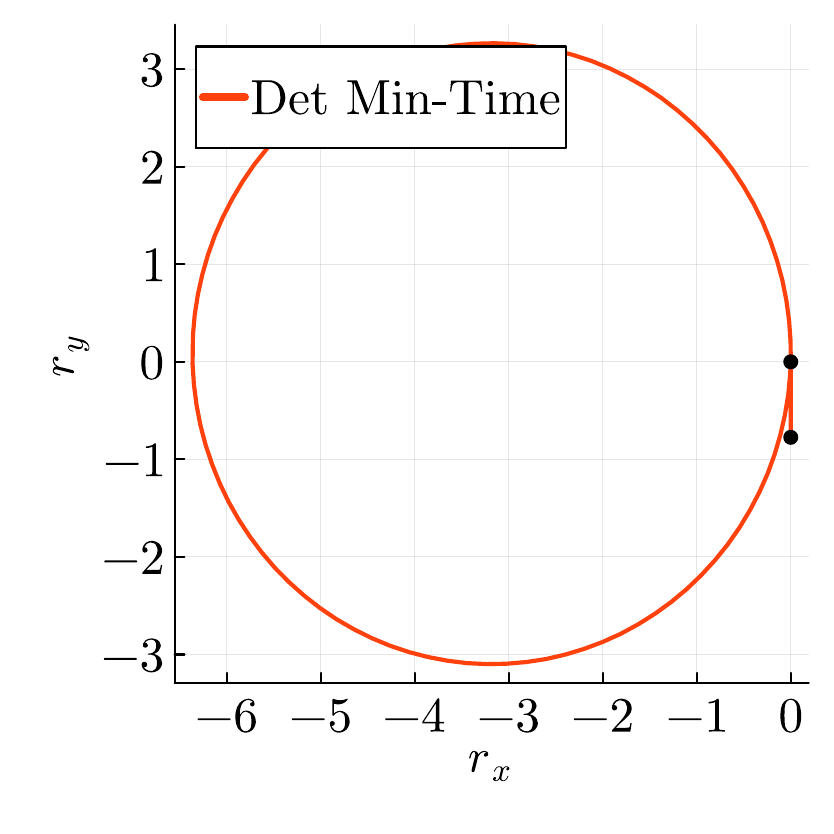}
         \label{fig:control}
     \end{subfigure}
    \caption{Plots of the trajectories using the proposed approach (\textit{Left}) and the deterministic minimum-time optimal control (\textit{Right}) in the representative scenario. The deterministic minimum-time optimal control passed $x_f$ once and reached $x_f$ through the arc of maximum curvature.}
       \label{fig:rep_dubin}
\end{figure}

As seen in Figure~\ref{fig:MCresults_dubin}, the mean terminal time-steps for the proposed approach is $34.7$, whereas that for the deterministic approach is $45.3$. Note that the distributions of terminal time-steps are largely divided into two groups, which are also analyzed in detail in Table~\ref{tab:monte_dubin}. The deterministic minimum-time optimal control generated a higher number of scenarios requiring more than $20$ steps to converge, which also increased the average of them. The primary reason for this bifurcation is that there are some scenarios where the vehicle overshoots the final destination $x_f$ within the consensus control horizon, resulting in a maneuver through the arc of maximum curvature to get back to the target. We clearly see that the deterministic minimum time-optimal control which doesn't consider initial state uncertainty results in such an overshooting more often than the proposed method. Figure~\ref{fig:rep_dubin} illustrates such a scenario, where the deterministic minimum-time optimal control takes longer times to converge than the proposed approach by taking this arc.

%% file: conclusion.tex
\section{Conclusion}
In this work, we presented an SCP-based framework to minimize the expected terminal time needed to reach a prescribed destination. This framework entails using a particle~\ac{MPC} approach with a consensus control horizon and a sum-of-norm cost. Numerical analyses for a spacecraft double integrator model and a Dubin's car model indicated that the proposed approach outperforms the deterministic approach successfully while being stable to parameters. Consequently, our proposed approach has two advantages: minimizing the expected terminal times under initial state distribution and efficient real-time computation of solutions offered by convex optimization-based techniques. One future approach could be to provide a sufficient condition of weight parameters in the sum-of-norm cost function such that the convexified optimization problem recovers the solution to the original minimum expected time optimal control problem. Finally, a future area of study includes incorporating our framework into a dual control approach that also seeks to minimize the initial state uncertainty within an optimal control framework~\cite {Feldbaum1963}.

%% file: main.bbl
\begin{thebibliography}{10}
\providecommand{\url}[1]{#1}
\csname url@rmstyle\endcsname
\providecommand{\newblock}{\relax}
\providecommand{\bibinfo}[2]{#2}
\providecommand\BIBentrySTDinterwordspacing{\spaceskip=0pt\relax}
\providecommand\BIBentryALTinterwordstretchfactor{4}
\providecommand\BIBentryALTinterwordspacing{\spaceskip=\fontdimen2\font plus
\BIBentryALTinterwordstretchfactor\fontdimen3\font minus \fontdimen4\font\relax}
\providecommand\BIBforeignlanguage[2]{{%
\expandafter\ifx\csname l@#1\endcsname\relax
\typeout{** WARNING: IEEEtran.bst: No hyphenation pattern has been}%
\typeout{** loaded for the language `#1'. Using the pattern for}%
\typeout{** the default language instead.}%
\else
\language=\csname l@#1\endcsname
\fi
#2}}

\bibitem{DeLaCroixRossiEtAl2024}
J.-P. {de} {la}~{Croix}, F.~Rossi, R.~Brockers, D.~Aguilar, K.~Albee, E.~Boroson, A.~Cauligi, J.~Delaune, and {others}, ``Multi-agent autonomy for space exploration on the {CADRE} {Lunar} technology demonstration mission,'' in \emph{IEEE Aerospace Conference}, 2024.

\bibitem{Harris20142}
M.~W. Harris and B.~A\c{c}\i{}kme\c{s}e, ``Minimum time rendezvous of multiple spacecraft using differential drag,'' \emph{{AIAA Journal of Guidance, Control, and Dynamics}}, vol.~37, no.~2, pp. 365--373, 2014.

\bibitem{Eren2017}
U.~Eren, A.~Prach, B.~B. Ko{\c{c}}er, S.~V. Rakovi{\'c}, E.~Kayacan, and B.~A{\c{c}}ikmese, ``Model predictive control in aerospace systems: Current state and opportunities,'' \emph{{AIAA Journal of Guidance, Control, and Dynamics}}, vol.~40, no.~7, pp. 1541--1566, 2017.

\bibitem{Richard2019}
R.~L. Sutherland, I.~V. Kolmanovsky, A.~R. Girard, F.~A. Leve, and C.~D. Petersen, ``Minimum-time model predictive spacecraft attitude controlfor waypoint following and exclusion zone avoidance,'' in \emph{{AIAA Scitech Forum}}, 2019.

\bibitem{CauligiSwanEtAl2023}
A.~Cauligi, R.~M. Swan, H.~Ono, S.~Daftry, J.~Elliott, L.~Matthies, and D.~Atha, ``{ShadowNav}: Crater-based localization for nighttime and permanently shadowed region lunar navigation,'' in \emph{{IEEE Aerospace Conference}}, 2023.

\bibitem{DaftryEtAl2023}
S.~Daftry, Z.~Chen, Y.~Cheng, S.~Tepsuporn, S.~Khattak, L.~Matthies, B.~Coltin, U.~Naam, L.~M. Ma, and M.~Deans, ``{LunarNav}: Crater-based localization for long-range autonomous rover navigation,'' in \emph{{IEEE Aerospace Conference}}, 2023.

\bibitem{Heath1987}
D.~Heath, S.~Orey, V.~Pestien, and W.~Sudderth, ``Minimizing or maximizing the expected time to reach zero,'' \emph{SIAM Journal on Control and Optimization}, vol.~25, no.~1, pp. 195--205, 1987.

\bibitem{Anderson2013}
R.~P. Anderson, E.~Bakolas, D.~Milutinovi\'{c}, and P.~Tsiotras, ``Optimal feedback guidance of a small aerial vehicle in a stochastic wind,'' \emph{Journal of Guidance, Control, and Dynamics}, vol.~36, no.~4, pp. 975--985, 2013.

\bibitem{GUILLOT2020148}
M.~Guillot and G.~Stauffer, ``The stochastic shortest path problem: A polyhedral combinatorics perspective,'' \emph{European Journal of Operational Research}, vol. 285, no.~1, pp. 148--158, 2020.

\bibitem{DyroHarrisonEtAl2021}
R.~Dyro, J.~Harrison, A.~Sharma, and M.~Pavone, ``Particle {MPC} for uncertain and learning-based control,'' in \emph{{IEEE/RSJ Int.\ Conf.\ on Intelligent Robots \& Systems}}, 2021.

\bibitem{BakoChenEtAl2011}
L.~Bako, D.~Chen, and S.~Lecoeuche. (2011) A numerical solution to the minimum-time control problem for linear discrete-time systems. {Available at }\url{https://arxiv.org/pdf/1109.3772.pdf}.

\bibitem{Malyuta2022}
D.~Malyuta, T.~P. Reynolds, M.~Szmuk, T.~Lew, R.~Bonalli, M.~Pavone, and B.~Açıkmeşe, ``Convex optimization for trajectory generation: A tutorial on generating dynamically feasible trajectories reliably and efficiently,'' \emph{IEEE Control Systems Magazine}, vol.~42, no.~5, pp. 40--113, 2022.

\bibitem{BlackmoreOnoEtAl2010}
L.~Blackmore, M.~Ono, A.~Bektassov, and B.~C. Williams, ``A probabilistic particle-control approximation of chance-constrained stochastic predictive control,'' \emph{{IEEE Transactions on Robotics}}, vol.~26, no.~3, pp. 502--517, 2010.

\bibitem{HOMEMDEMELLO2014}
T.~H. de~Mello and G.~Bayraksan, ``{Monte} {Carlo} sampling-based methods for stochastic optimization,'' \emph{Surveys in Operations Research and Management Science}, vol.~19, no.~1, pp. 56--85, 2014.

\bibitem{ZhuAlonsoMora2019}
H.~Zhu and J.~Alonso-Mora, ``Chance-constrained collision avoidance for {MAVs} in dynamic environments,'' \emph{{IEEE Robotics and Automation Letters}}, vol.~4, no.~2, pp. 776--783, 2019.

\bibitem{Weissel2009}
F.~Weissel, M.~F. Huber, and U.~D. Hanebeck, \emph{Stochastic Nonlinear Model Predictive Control based on Gaussian Mixture Approximations}.\hskip 1em plus 0.5em minus 0.4em\relax Berlin, Heidelberg: Springer Berlin Heidelberg, 2009, pp. 239--252.

\bibitem{SinghSingla2007}
T.~Singh and P.~Singla, ``Sequential linear programming for design of time-optimal controllers,'' in \emph{{Proc.\ IEEE Conf.\ on Decision and Control}}, 2007.

\bibitem{candes2008}
E.~J. Cand{\`e}s, M.~B. Wakin, and S.~P. Boyd, ``Enhancing sparsity by reweighted {$\mathcal{L}$}1 minimization,'' \emph{Journal of Fourier Analysis and Applications}, vol.~14, no.~5, pp. 877--905, 2008.

\bibitem{Nagahara2016}
M.~Nagahara, D.~E. Quevedo, and D.~Nešić, ``Maximum hands-off control: A paradigm of control effort minimization,'' \emph{IEEE Transactions on Automatic Control}, vol.~61, no.~3, pp. 735--747, 2016.

\bibitem{Kazu2023LCSS}
K.~Echigo, C.~R. Hayner, A.~Mittal, S.~B. Sarsılmaz, M.~W. Harris, and B.~Açıkmeşe, ``Linear programming approach to relative-orbit control with element-wise quantized control,'' \emph{IEEE Control Systems Letters}, vol.~7, pp. 3042--3047, 2023.

\bibitem{DomahidiChuEtAl2013}
A.~Domahidi, E.~Chu, and S.~Boyd, ``{ECOS}: An {SOCP} solver for embedded systems,'' in \emph{{European Control Conference}}, 2013.

\bibitem{Natarajan1995}
B.~K. Natarajan, ``Sparse approximate solutions to linear systems,'' \emph{SIAM Journal on Computing}, vol.~24, no.~2, pp. 227--234, 1995.

\bibitem{Candes2005}
E.~Candes and T.~Tao, ``Decoding by linear programming,'' \emph{IEEE Transactions on Information Theory}, vol.~51, no.~12, pp. 4203--4215, 2005.

\bibitem{ChariKamathEtAl2024}
G.~M. Chari, A.~G. Kamath, P.~Elango, and B.~A\c{c}{\i}kme\c{s}e, ``Fast {Monte} {Carlo} analysis for 6-{DoF} powered-descent guidance via {GPU}-accelerated sequential convex programming,'' in \emph{{AIAA Scitech Forum}}, 2024.

\bibitem{MoteEgerstedtEtAl2020}
M.~Mote, M.~Egerstedt, E.~Feron, A.~Bylard, and M.~Pavone, ``Collision-inclusive trajectory optimization for free-flying spacecraft,'' \emph{{AIAA Journal of Guidance, Control, and Dynamics}}, 2020.

\bibitem{CauligiChakrabartyEtAl2022}
A.~Cauligi, A.~Chakrabarty, S.~{Di}~{Cairano}, and R.~Quirynen, ``{PRISM}: Recurrent neural networks and presolve methods for fast mixed-integer optimal control,'' in \emph{{Learning for Dynamics \& Control}}, 2022.

\bibitem{Dubin1957}
L.~E. Dubins, ``On curves of minimal length with a constraint on average curvature, and with prescribed initial and terminal positions and tangents,'' \emph{American Journal of Mathematics}, vol.~79, no.~3, pp. 497--516, 1957.

\bibitem{taylor2021}
T.~P. Reynolds, ``Computation guidance and control for aerospace systems,'' Ph.D. dissertation, University of Washington, Seattle, WA, 2021.

\bibitem{szmuk2019}
M.~Szmuk, ``Successive convexification \& high performance feedback control for agile flight,'' Ph.D. dissertation, University of Washington, Seattle, WA, 2019.

\bibitem{Feldbaum1963}
A.~A. Feldb{\^a}um, ``Dual control theory problems,'' \emph{IFAC Proceedings Volumes}, vol.~1, no.~2, pp. 541--550, 1963.

\end{thebibliography}
